# Combined Heat and Power Unit Commitment with Smart Parking Lots of Plug-in Electric Vehicles

Hamidreza Sadeghian, Zhifang Wang
Department of Electrical and Computer Engineering
Virginia Commonwealth University, Richmond, VA, USA
Email: {sadeghianh}, {zfwang} @vcu.edu

*Abstract*—**Vehicle-to-grid (V2G) technology has drawn great interest in the recent years and its efficiency depends on scheduling of charging process of plug-in electric vehicles (PEVs) as small portable power plants in smart parking lots. On the other hand, active shift from centralized to decentralized power generation and environmental concerns have caused an increase in the utilization of combined heat and power (CHP) units in power systems. The goal of this study is to develop and simulate a novel approach for combined heat and power unit commitment with PEVs for cost reduction in electric power system. A schedule for charging and discharging processes of PEVs with respect to load curve variations is proposed in this paper. A modified test system consisting of conventional TG (thermal power generating) units, CHP units, and PEVs is employed to investigate the impacts of PEVs on generation scheduling. For the modified test system, results obtained are encouraging and indicate both the feasibility of the proposed technique and its effectiveness on generation scheduling.**

## I. INTRODUCTION

Smart grid incorporates a communication infrastructure that enables system components to exchange information and commands securely and reliably. By providing such real time information, independent system operators (ISOs) achieve the ability to manage generating resources associated with the load demand concurrently. Online supervision may cause to handle supply demand equilibrium in a real time fashion and better integration of renewable energies associated with electricity storages [1, 2]. PEVs as portable electricity storages are potentially not only environmentally friendly and quiet but also cost-effective in terms of operating costs and energy prices compared to conventional TG units. PEVs can decrease dependencies on small expensive units by discharging at peak hours through participating in vehicle to grid (V2G) service which leads to operating costs reduction [3]. Moreover other applications of PEVs include flattening load curve and maximizing loadability of the network by discharging PEVs in peak time and charging during off-peak periods [4]. Other than reducing emission and mitigating transmission line congestion, PEVs have ability for fast load tracking due to high availability of parking lots and local management of generation subject to placing parking lots near the load centers [5].

In addition, to reduce dependency on external energy supplies and mitigate climate change, many countries have adopted policies to increase both energy conservation and the share of renewable energy resources [6]. The benefit of renewable energy aggregation has been demonstrated in a number of studies [7]. Moreover, in some countries and regions e.g. in the EU such policies include increasing the share of combined heat and power (CHP) [8]. CHP or cogeneration is the production of electricity and thermal energy in a single, integrated system. It is a valuable energy production technology that can yield much higher total energy efficiency than separate heat and power generation. The fuel efficiency of CHP production unit can be as much as 90%. In fact high efficiency of CHP systems is one of the major factors that make them attractive for investors [9].

Several studies are reported on impacts of PEVs on generation scheduling problem. Authors in [10] developed a unit commitment model for Texas electric power system considering plug-in hybrid electric vehicle fleet to measure the potential cost savings. An intelligent unit commitment associated with vehicle to grid has been studied based upon minimizing operating costs in addition to environmental emission externalities [11]. In [12], a charging and discharging schedule of PEVs with respect to load curve variations with thermal generation scheduling is proposed. Literature review reveals a gap for studying the impacts of CHP systems and PEVs integration in conventional TG units and scheduling.

In this paper, a bridge between PEVs and CHP unit commitment problem has been made and the impacts of parking lots penetration on combined heat and power generation scheduling in smart grids have been investigated. The contributions of this study can be summarized as follows: 1) A new structure of combined heat and power unit commitment incorporating PEVs, so-called CHPUC-PEV was proposed. 2) A modified double Benders decomposition method was utilized to solve the proposed optimization problem. 3) The performance and effectiveness of the proposed approach are evaluated with numerical simulations.

## II. CHPUC-PEV FORMULATION

### A. Nomenclature

| | |
|---|---|
| $F(P_{i,t}, H_{i,t})$ | fuel consumption function of $i$th unit at time $t$ |
| $HD_t$ | total system heat demand at time $t$ |
| $H_i^{max}$ | maximum heat output of $i$th generator (MW) |
| $H_i^{min}$ | minimum heat output of *i*th generator (MW) |
| $H_{i,t}$ | heat level of *i*th generator at *t*th hour (MW) |
| $i$ | index for generator unit |
| $j$ | index for smart parking lot |
| $M$ | total number of smart parking lots |
| $N$ | total number of all generating units |
| $N_{cov}$ | total number of conventional TG units |
| $N_{heat}$ | total number of heat only units |
| $N_{dsch,t}^j$ | # of connected discharging vehicles to the grid at hour t at smart parking lot *j* |
| $N_{ch,t}^j$ | # of connected charging vehicles to the grid at hour t at smart parking lot *j* |
| $N^{max}$ | total vehicles in the system |
| $N_{dsch}^{max}$ | minimum # of discharging vehicles at hour t |
| $N_{dsch}^{min}$ | maximum # of discharging vehicles at hour t |
| $N_{ch}^{min}$ | minimum # of charging vehicles at hour t |
| $N_{ch}^{max}$ | maximum # of charging vehicles at hour t |
| $PD_t$ | total system power demand at time *t* |
| $P_i^{max}$ | maximum power output of *i*th unit (MW) |
| $P_i^{min}$ | minimum power output of *i*th unit (MW) |
| $P_{i,t}$ | power level of *i*th generator at *t*th hour (MW) |
| $\pi_j$ | V2G cost coefficient for smart parking lot *j* |
| $P_{PEV,j}$ | available PEV power for V2G in parking lot *j* |
| $pv_j$ | capacity of each vehicle for smart parking lot *j* |
| $RD_t$ | total system reserve demand at time *t* |
| $SD_{i,t}$ | shutdown costs of thermal unit *i* at time *t* |
| $SH_i$ | shutdown cost of unit *i*($) |
| $ST_i$ | startup cost of unit *i*($) |
| $SU_{i,t}$ | startup costs of thermal unit *i* at time *t* |
| $T$ | dispatch period in hours |
| $t$ | index for time |
| $T_i^{down}$ | minimum down-time of *i*th generator unit |
| $t_{i,t}^{down}$ | time duration for *i*th unit that has been OFF at *t* |
| $T_i^{up}$ | minimum up-time of *i*th generator unit |
| $t_{i,t}^{up}$ | time duration for ith unit that has been ON at *t* |
| $X_{i,t}$ | commitment state of *i*th unit at *t*th hour |
| $\delta_j$ | state of charge for smart parking lot *j* |
| $\eta_j$ | total efficiency for smart parking lot *j* |

*B. Objective function*

The objective of CHPUC-PEV problem is minimizing the total operation cost of units over a scheduling period. Indeed, the CHPUC-PEV objective function usually includes different terms such as the start-up cost, shutdown cost and the fuel cost, while, the term of PEVs energy costs has been also taken into account in this manuscript.

$$\min_{X,P,H,N_{dsch}^j} f = \sum_{i=1}^{N}\sum_{t=1}^{T} F_i(P_{i,t}, H_{i,t}).X_{i,t} + SU_{i,t} + SD_{i,t} + \sum_{j=1}^{M}\sum_{t=1}^{T} N_{dsch,t}^j . P_{PEV,j} . \pi_j \quad (1)$$

where $SU_{i,t}$ and $SD_{i,t}$ represent startup and shutdown costs of thermal unit *i* at time *t*, respectively, which are determined based on the following inequalities:

$$SU_{i,t} = ST_i \times X_{i,t} \times (X_{i,t} - X_{i,t-1}) \quad (2)$$
$$SD_{i,t} = SH_i \times X_{i,t-1} \times (X_{i,t-1} - X_{i,t}) \quad (3)$$

In Eq. (1), second part represents PEVs usage cost as small portable power plants. In this term, $P_{PEV,j}$ defines as available PEV power for V2G in smart parking lot *j* which can be given by the Eq. (2).

$$P_{PEV,j} = pv_j . \delta_j . \eta_j \quad (4)$$

where $pv_j$, $\eta_j$, $\delta_j$, represent the average capacity of each battery, inverter efficiency parameter, and departure state of charge parameter for smart parking lot *j*, respectively. In this paper we considered hourly scheduling of unit commitment and PEV charging/discharging. Therefore, the complexities related with fast dynamics are ignored and a simple statistic model is used instead which is sufficient for our study. The objective function includes generation cost of both thermal units and CHP units. The fuel costs are presented as follows:

Thermal units

$$F_i(P_{i,t}) = a_i + b_i.P_{i,t} + c_i.(P_{i,t})^2 \quad (5)$$

CHP units

$$F_i(P_{i,t}) = a_i + b_i.P_{i,t} + c_i.(P_{i,t})^2 + d_i.H_{i,t} + e_i.(H_{i,t})^2 + f_i.P_{i,t}.H_{i,t} \quad (6)$$

where $a_i$, $b_i$, $c_i$, $d_i$, $e_i$ and $f_i$ are positive fuel cost coefficients of *i*th unit.

Fig. 1 shows the heat-power feasible operation region (FOR) of a CHP unit. It should be mentioned that CHP units assumed to utilize water vapor or gas turbines and both of them are modeled by means of a feasible operation region designated by ABCDEF that links the electric power generation and heat production [13]. Along the boundary curve BC, the heat capacity increases as the power generation decrease. Although in most cases the fuel cost is a convex function, the feasible operation region of advanced CHP units is non-convex. From Fig. 1, in CHP systems power generation depends on the heat generation and vice versa, which implies that the production of heat and power must be planned in coordination with each other.

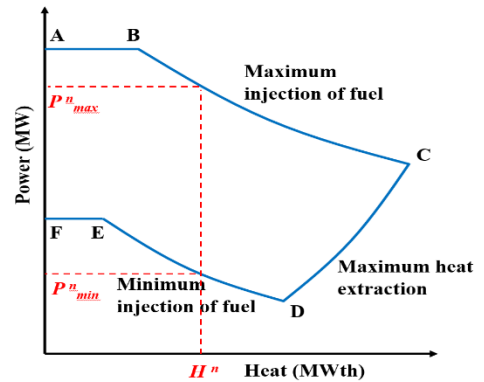

Figure 1. Feasible operation region for a CHP unit

## C. Constraints

The CHPUC-PEV involves constraints like power and heat balance, spinning reserve requirement, minimum up/down time of a unit and limited number of PEVs. In addition, a unit is to generate power within a given range. The formulations of the constraints are given below:

**Power and heat balance:** The power and heat generated by all the committed units with PEVs power at a time instant must meet power and heat demand at that time instant, respectively.

$$\sum_{i=1}^{N} P_{i,t} X_{i,t} + \sum_{j=1}^{M} (pv.N_{dsch,t}^{j} - pv.N_{ch,t}^{j}) = PD_t \quad (7)$$

$$\sum_{i=1}^{N} H_{i,t} X_{i,t} = HD_t \quad (8)$$

**Generation limit:** The generation limit of a CHP unit is specified by the FOR as shown in Section II.B and Fig.1. The upper and lower limits of the conventional TG units and the heat only units can be described by constraints (9) and (10) respectively as follows:

$$P_i^{\min} \leq P_{i,t} \leq P_i^{\max}, \; i=1,\cdots,N_{cov} \quad (9)$$

$$H_j^{\min} \leq H_{j,t} \leq H_j^{\max}, \; j=1,\cdots,N_{heat} \quad (10)$$

**Minimum up and down times:** It is considered that a unit must be on/off for a minimum time before it can be shut down or restarted, respectively:

$$(t_{i,t}^{up} - T_i^{up})(X_{i,t-1} - X_{i,t}) \geq 0 \quad (11)$$

$$(t_{i,t}^{down} - T_i^{down})(X_{i,t} - X_{i,t-1}) \geq 0 \quad (12)$$

**Spinning reserve:** In order to fast response to compensate the deviation between real and predictive demand in power system, spinning reserve is required. Mathematically, spinning reserve requirement at each hour is the total amount of maximum capacity of all synchronized units minus the total generating output in that hour which can be given by the Eq. (13).

$$\sum_{i=1}^{N} P_{i,\max} X_{i,t} + \sum_{j=1}^{M} (pv.N_{dsch,t}^{j} - pv.N_{ch,t}^{j}) \geq PD_t + RD_t \quad (13)$$

**Charge/discharge limits:** In supporting the daily use and have a reliable operation, a certain amount of power should be injected into the PEVs batteries. Therefore, it is necessary to limit the accumulated charging/discharging power. In addition, limited number of PEVs should charge/discharge at the same time over a predefined horizon.

$$\sum_{j=1}^{M} \sum_{t=1}^{T} N_{dsch,t}^{j} = N^{\max} \quad (14)$$

$$N_{dsch}^{\min} \leq N_{dsch,t}^{j} \leq N_{dsch}^{\max} \quad (15)$$

$$\sum_{j=1}^{M} \sum_{t=1}^{T} N_{ch,t}^{j} = N^{\max} \quad (16)$$

$$N_{ch}^{\min} \leq N_{ch,t}^{j} \leq N_{ch}^{\max} \quad (17)$$

In this study, charging/discharging frequency is assumed once a day. Each vehicle should have a desired departure state of charge (SOC) level, while $\eta$ is defined as integrated efficiency for charging/discharging plus inverter.

## III. DOBBLE BENDERS DECOMPOSITION APPROACH

As a mixed-integer nonlinear optimization problem with a large number of continuous and integer control variables, the CHPUC-PEV problem is non-convex due to CHP feasible operation region. In this study, double Benders decomposition (DBD) approach is used to solve the proposed CHPUC-PEV problem. In general, Benders decomposition (BD) approach decomposes the original problem into one master problem and several sub-problems. By solving each sub-problem, a set of dual variables are obtained and used to generate benders cuts for the master problem. The procedure to implement the proposed DBD approach for solving CHPUC-PEV problem is shown in Fig. 2, where it is consist of two BD algorithms, namely the outer BD and the inner BD. For the outer BD, the master problem determines the integer variables (on/off state of each generating unit) and the sub-problem solves the economic dispatch (ED) along with charge/discharge scheduling problem.

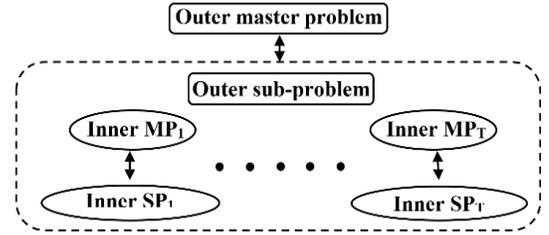

Figure 2. Strategy of the proposed double Benders decomposition approch

As it mentioned, for most CHP units, the heat production capacities depend on the power generation and vice versa, which implies that the production of heat and power must be planned in coordination with each other. Therefore, it presents a natural decomposition scheme for the BD algorithm. For most CHP systems feasible operation region is non-convex and it increases the complexity to get the global optimal solution. However, decomposing heat and power variables with BD approach results to a simple convexification procedure and then we solve convex problems in both master problem and sub-problem which it leads to find global optimum solution [14, 15]. In the proposed approach in this study, for solving the ED problem (outer sub-problem) another BD algorithm (inner BD) is used. For this BD algorithm, the variables representing the heat production are solved for in the master problem while the participating PEVs in addition to the ones representing the power production are kept in the sub problem. More details of the proposed DBD approach can be found in [14].

The objective function is decomposed to the objective functions of the master problem and sub-problem. The master problem in the proposed solution approach has mixed integer linear programming model with the following objective function

$$\min_{X} \; \mu_M = \left\{ \sum_{t=1}^{T} \left[ \sum_{i=1}^{N} \left( (a_i.X_{i,t}) + SD_{i,t} + SU_{i,t} \right) + \gamma_t \right] \right\} \quad (18)$$

where $\gamma_t$ denotes the benders cuts obtained from dual variables and included constraints for the outer master problem are the inequalities given by Eqs. (2) and (3) as well as Eqs. (11) and (12).

In the outer sub-problem, the dispatch decision variables, $P_{i,t}$, $H_{i,t}$, and $N_{dsch,t}$ are accounted. An ED allocates power and heat generation among committed units based on a BD algorithm. The objective function of inner sub-problem is as follows

$$\min_{P, N_{dsch}} RF(X_t, H_t, P_t, N_{dsch,t}^j)$$
subject to (19)
$$X_t = X_t^{(\kappa)} : \lambda x_t^{(\kappa)}$$
$$H_t = H_t^{(v)} : \lambda h_t^{(v)}$$

where $H_t^{(v)}$ represents the value of heat production vector at $v^{th}$ iteration of inner BD, $X_t^{(\kappa)}$ is the value of units status vector at $\kappa^{th}$ iteration of outer BD and $\lambda$ is prefix indicating dual variables. The additional constraints of the inner sub-problem include the power demand inequality (Eq. (7)), power generation limit (Eqs. (18)), minimum spinning reserve constraint (Eq. (13)), and charge/discharging limits (Eqs. (14-17)). The output of this sub-problem determines the values of participating PEV, $P_t^{(v)}$ and the dual variable vector associated with those constraints that fix the complicating variables, $H_t$. The objective function of inner master problem solely as a function of the complicating variables ($H$) for time step $t$ is given by

$$\min \quad \mu_{MM,t} \quad (20)$$

subject to

$$\mu_{MM,t} \geq RF(X_t^{(\kappa)}, P_t^{(\sigma)}, H_t^{(\sigma)}, N_{dsch,t}^j) + \sum_{i=1}^{N} \lambda h_{i,t}^{(\sigma)} \left( H_{i,t} - H_{i,t}^{(\sigma)} \right) \quad (21)$$

where σ indicates iteration with the highest value of normal Benders cut. Assumed constraints of the inner master problem include the heat demand balance (Eq. (8)), and heat generation limits (Eq. (10)). In this study, strong Benders cut is used which can enhance the convergence of the BD approach proposed [16]. If a sub-problem becomes infeasible, the strong Benders cut is added to the master problem of the next iteration as well, however, despite of the previous case (feasible sub-problem), comparison among iterations does not include the current iteration. Solving ED problem for each hour provides Benders cuts ($\gamma_t$) for outer master problem,

$$\gamma_t \geq RF(X_t^{(\kappa)}, P_t^{(v)}, H_t^{(v)}, N_{dsch,t}^{j(v)}) + \sum_{i=1}^{N} \lambda x_{i,t}^{(\kappa)} \left( X_{i,t} - X_{i,t}^{(\kappa)} \right) \quad (22)$$

## IV. SIMULATION AND RESULTS

The modified 11 unit system used in the simulation of PEVs impact on generation scheduling is based on data presented in [14]. A total number of 50000 PEVs aggregated from multiple smart parking lots in the system are considered the simulations. Spinning reserve requirement is assumed to be 10% of the hourly load demand in 24 hour scheduling time period. Following parameters are assumed for PEVs: maximum battery capacity = 25 kWh, average battery capacity = 15 kWh, minimum battery capacity = 10 kWh, charging/discharging frequency = 1 per day, total efficiency (η) = 85%, and state of charge (δ) = 50%.

In this paper, three different scenarios are investigated (Table I). First scenario consists of a typical cost based unit commitment problem with obligation of heat and power demand. In the second scenario, an integration of 50,000 PEVs charged by renewable sources as the CHPUC-PEV problem is considered. It is assumed that 50,000 PEVs are aggregated in parking lots managed by an aggregator organization and supplied by renewable energies such as wind turbines and solar power, near the parking lot to avoid transmission losses. The third scenario studies the integration of 50,000 PEVs charged/discharged via power grid without renewable sources in the CHPUC-PEV problem.

It is worthy to mention that the "charged" state is when the PEVs are charged via the electrical grid and the "discharged" state is when batteries of PEVs are depleted to deliver their stored energy to the electrical grid. Maximum number of charging/discharging vehicles at each hour for scenarios 2 and 3 is 10% and 20% of total vehicles, respectively.

TABLE I. SUMMARY OF SCENARIOS# 1-3.

| Scenarios | Scenario definition |
|---|---|
| Scenario#1 | Thermal and CHP units without PEVs |
| Scenario#2 | Thermal and CHP units with PEVs and renewable energy |
| Scenario#3 | Thermal and CHP units with PEVs and without renewable energy |

### A. Scenario 1

In this scenario, a modified 11 unit system consists of eight conventional TG units, two CHP units and one heat-only unit is utilized. For 24 hours of a day, the total cost obtained $914,008.22 as shown in Table II. The best results for the three scenarios are compared with each other in this table. Note that, the large difference between the best cost of the third scenario and the second scenario is the result of installing additional renewable generation units in in the second scenario.

TABLE II. RESULTS OF PROPOSED METHOD FOR DIFFERENT SCENARIOS

| Scenarios | Best result ($) | Improvement ($) | Improvement (%) |
|---|---|---|---|
| Scenario#1 | 914,008.22 | - | - |
| Scenario#2 | 904,819.22 | 9,189.00 | 1.00 |
| Scenario#3 | 911,488.60 | 2,519.62 | 0.27 |

### B. Scenario 2

In this case, we considered that 10-units (Thermal and CHP units) would be employed to supply the grid electricity and in addition, CHP units plus one heat-only unit will meet the heat load demands for the test system. In addition, PEVs are totally supplied by renewable sources. In this regard, the best operating cost equals to $904,819.22. By discharging PEVs during particular hours and supplying adequate power, operation hours of committed units are decreased and consequently total operation cost is improved. Moreover, minimum number of discharging PEVs for each hour is assumed as shown in Table III. The best results of the second scenario is presented in Table II to enable comparison with other two scenarios. It is worthwhile to note that in this study one smart parking lot is

considered, however, by considering multiple penetration of PEVs parking lots, higher improvements will be achievable.

TABLE III. MINIMUM NUMBER OF DISCHARGING PEVS.

| Hour | $N_{dsch}^{min}$ | Hour | $N_{dsch}^{min}$ | Hour | $N_{dsch}^{min}$ | Hour | $N_{dsch}^{min}$ |
|---|---|---|---|---|---|---|---|
| 1 | 0 | 7 | 0 | 13 | 3400 | 19 | 0 |
| 2 | 0 | 8 | 0 | 14 | 1500 | 20 | 0 |
| 3 | 2000 | 9 | 1500 | 15 | 0 | 21 | 0 |
| 4 | 0 | 10 | 3400 | 16 | 0 | 22 | 0 |
| 5 | 2200 | 11 | 3400 | 17 | 0 | 23 | 0 |
| 6 | 0 | 12 | 3400 | 18 | 0 | 24 | 0 |

*C. Scenario 3*

Figure. 3 illustrates the heat and power load curve of the test system, where it is divided into three different periods based on the power load, i.e. valley period, off-peak period, and peak period . To reduce the operation costs, PEVs are charged during off-peak periods and discharged during peak periods. The result of unit commitment problem in presence of PEVs are presented in Tables IV and V in the appendix. Table IV shows the best result of the second scenario and the best results of the third scenario is shown in Table V. As it can be seen from the results, costs of all scenarios that include PEVs in unit commitment problem are lower than conventional CHP unit commitment which shows the effectiveness of incorporating PEVs in total operation cost reduction. It is worthwhile to note that

, the main challenge of unit commitment is to properly schedule small expensive units, as large cheap units are always on. Operators expect that large cheap units will mainly satisfy base load and other small expensive units will fulfill the peak loads. Using PEVs as small portable power plants reduce dependencies on small expensive units.

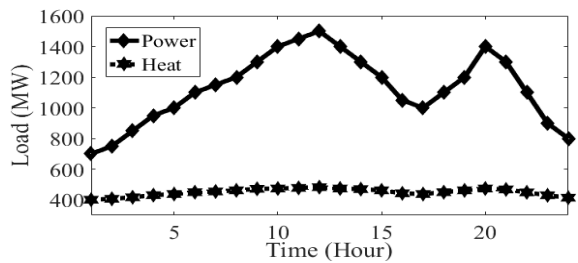

Figure 3. Heat and power load curve of the test system [12]

## V. CONCLUSION

This paper presents a new approach to solve combined heat and power unit commitment problem in presence of plug-in electric vehicles as small portable generation units in addition to typical generation constraints in unit commitment problem. The numerical results illustrate the effectiveness of this approach on the system operation cost reduction. In fact, PEVs not only eliminate the need for small expensive units in the power systems, but also provide additional reserve capacity and reliability for the existing power systems. The results of generation scheduling with the optimal charging/discharging scheme of PEVs suggests that the total operation cost decreases noticeably without any change in the total energy consumption. However, future work will incorporate multiple penetration of PEVs smart parking lots with different settings of PEVs operations. In addition, the integration of uncertain renewable energy generation is also worth investigation.

APPENDIX

TABLE IV. SCHEDULE AND DISPATCH OF GENERATING UNITS FOR SCENARIO 2 (TOTAL OPERATION COST = **$904,819.22**)

| Hour | Power (MW) | | | | | | | | | | Heat (MWth) | | | $N_{dsch}$ | $v2g$ $Pv2g$ (MW) | Hourly cost ($) |
|---|---|---|---|---|---|---|---|---|---|---|---|---|---|---|---|---|
| | U1 | U2 | U3 | U4 | U5 | U6 | U7 | U8 | U9 | U10 | U9 | U10 | Boiler | | | |
| 1 | 455 | 115.77 | 0 | 0 | 0 | 0 | 0 | 0 | 81 | 40 | 105 | 75 | 219 | 1291 | 8.23 | 28,450 |
| 2 | 455 | 171.6 | 0 | 0 | 0 | 0 | 0 | 0 | 81 | 40 | 105 | 75 | 206 | 376 | 2.40 | 29,182 |
| 3 | 455 | 257.22 | 0 | 0 | 0 | 0 | 0 | 0 | 82 | 43 | 104 | 75 | 238 | 2006 | 12.78 | 31,449 |
| 4 | 455 | 351.62 | 0 | 0 | 0 | 0 | 0 | 0 | 81 | 46 | 105 | 75 | 279 | 2569 | 16.38 | 34,054 |
| 5 | 455 | 398.76 | 0 | 0 | 0 | 0 | 0 | 0 | 82 | 47 | 104 | 75 | 287 | 2704 | 17.25 | 35,153 |
| 6 | 455 | 455 | 0 | 60.76 | 0 | 0 | 0 | 0 | 81 | 44 | 105 | 75 | 298 | 665 | 4.26 | 37,276 |
| 7 | 455 | 455 | 0 | 98.83 | 0 | 0 | 0 | 0 | 82 | 47 | 104 | 75 | 304 | 3164 | 20.18 | 38,095 |
| 8 | 455 | 455 | 130 | 42.75 | 0 | 0 | 0 | 0 | 82 | 40 | 104 | 75 | 311 | 118 | 0.76 | 39,255 |
| 9 | 455 | 455 | 130 | 120.26 | 0 | 0 | 0 | 0 | 81 | 48 | 105 | 75 | 320 | 1685 | 10.75 | 41,169 |
| 10 | 455 | 455 | 130 | 130 | 67.67 | 0 | 0 | 0 | 81 | 52 | 105 | 75 | 322 | 4604 | 29.36 | 43,025 |
| 11 | 455 | 455 | 130 | 130 | 77.46 | 40 | 0 | 0 | 81 | 53 | 105 | 88 | 313 | 4477 | 28.55 | 44,241 |
| 12 | 455 | 455 | 130 | 130 | 75.54 | 50 | 44.23 | 0 | 81 | 55 | 105 | 88 | 318 | 3801 | 24.24 | 45,883 |
| 13 | 455 | 455 | 130 | 130 | 58.68 | 40 | 0 | 0 | 82 | 46 | 105 | 79 | 318 | 4499 | 28.69 | 43,138 |
| 14 | 455 | 455 | 130 | 123.89 | 0 | 0 | 0 | 0 | 81 | 45 | 105 | 78 | 314 | 1586 | 10.11 | 40,967 |
| 15 | 455 | 455 | 130 | 35.57 | 0 | 0 | 0 | 0 | 81 | 42 | 105 | 75 | 310 | 224 | 1.43 | 39,252 |
| 16 | 455 | 455 | 15 | 0 | 0 | 0 | 0 | 0 | 81 | 43 | 105 | 75 | 291 | 157 | 1.00 | 36,186 |
| 17 | 455 | 455 | 53.89 | 0 | 0 | 0 | 0 | 0 | 81 | 46 | 105 | 75 | 288 | 2649 | 16.82 | 35,035 |
| 18 | 455 | 455 | 38.89 | 0 | 0 | 0 | 0 | 0 | 82 | 49 | 104 | 75 | 299 | 3155 | 20.11 | 37,084 |
| 19 | 455 | 455 | 130 | 30.75 | 0 | 0 | 0 | 0 | 82 | 44 | 104 | 75 | 311 | 510 | 3.25 | 39,311 |
| 20 | 455 | 455 | 130 | 130 | 55 | 44.75 | 0 | 0 | 82 | 48 | 105 | 75 | 326 | 39 | 0.24 | 43,965 |
| 21 | 455 | 455 | 130 | 62.54 | 45.9 | 0 | 0 | 0 | 82 | 47 | 104 | 75 | 317 | 3539 | 22.57 | 41,156 |
| 22 | 455 | 455 | 0 | 0 | 50.76 | 0 | 0 | 0 | 81 | 43 | 105 | 75 | 297 | 2391 | 15.24 | 37,321 |
| 23 | 455 | 311.34 | 0 | 0 | 0 | 0 | 0 | 0 | 81 | 44 | 105 | 75 | 278 | 1358 | 8.66 | 33,246 |
| 24 | 455 | 208.13 | 0 | 0 | 0 | 0 | 0 | 0 | 81 | 40 | 105 | 75 | 262 | 2489 | 15.87 | 30,926 |

TABLE V. SCHEDULE AND DISPATCH OF GENERATING UNITS FOR SCENARIO 3 (TOTAL OPERATION COST = **$911,488.60**)

| Hour | Power (MW) | | | | | | | | | | Heat (MWth) | | | $N_{dsch}$ | $Pv2g$ (MW) | $N_{ch}$ | $Pg2v$ (MW) | Hourly cost ($) |
|---|---|---|---|---|---|---|---|---|---|---|---|---|---|---|---|---|---|---|
| | U1 | U2 | U3 | U4 | U5 | U6 | U7 | U8 | U9 | U10 | U9 | U10 | Boiler | | | | | |
| 1 | 455 | 153.11 | 0 | 0 | 0 | 0 | 0 | 0 | 81 | 40 | 105 | 75 | 220 | 0 | 0 | 4566 | 29.10 | 29097 |
| 2 | 455 | 219.72 | 0 | 0 | 0 | 0 | 0 | 0 | 81 | 40 | 105 | 75 | 206 | 0 | 0 | 7172 | 45.81 | 30019 |
| 3 | 455 | 318.61 | 0 | 0 | 0 | 0 | 0 | 0 | 82 | 43 | 104 | 75 | 238 | 0 | 0 | 7625 | 48.63 | 32513 |
| 4 | 455 | 379.28 | 0 | 0 | 0 | 0 | 0 | 0 | 81 | 50 | 105 | 75 | 279 | 0 | 0 | 2397 | 15.30 | 34702 |
| 5 | 455 | 419.48 | 0 | 0 | 0 | 0 | 0 | 0 | 82 | 52 | 104 | 75 | 287 | 0 | 0 | 1330 | 8.80 | 35722 |
| 6 | 455 | 455 | 0 | 85.37 | 0 | 0 | 0 | 0 | 81 | 51 | 105 | 75 | 298 | 0 | 0 | 4290 | 27.36 | 37993 |
| 7 | 455 | 455 | 0 | 118.65 | 0 | 0 | 0 | 0 | 82 | 50 | 104 | 75 | 304 | 0 | 0 | 1671 | 10.24 | 38678 |
| 8 | 455 | 455 | 130 | 67.82 | 0 | 0 | 0 | 0 | 82 | 40 | 104 | 75 | 311 | 0 | 0 | 4678 | 29.80 | 39791 |
| 9 | 455 | 455 | 130 | 78.72 | 0 | 0 | 0 | 0 | 81 | 53 | 105 | 75 | 320 | 7416 | 47.29 | 0 | 0 | 40672 |
| 10 | 455 | 455 | 130 | 130 | 41.6 | 0 | 0 | 0 | 81 | 52 | 105 | 75 | 322 | 8690 | 55.41 | 0 | 0 | 42546 |
| 11 | 455 | 455 | 130 | 130 | 71.7 | 40 | 25.96 | 0 | 81 | 53 | 105 | 88 | 313 | 1308 | 8.34 | 0 | 0 | 44816 |
| 12 | 455 | 455 | 130 | 130 | 73.21 | 50 | 48.93 | 10 | 81 | 55 | 105 | 88 | 318 | 1860 | 11.86 | 0 | 0 | 46236 |
| 13 | 455 | 455 | 130 | 130 | 46 | 0 | 0 | 0 | 82 | 52 | 105 | 79 | 318 | 7550 | 48.14 | 0 | 0 | 42939 |
| 14 | 455 | 455 | 130 | 124.75 | 0 | 0 | 0 | 0 | 81 | 53 | 105 | 78 | 314 | 196 | 1.25 | 0 | 0 | 41310 |
| 15 | 455 | 455 | 130 | 52.59 | 0 | 0 | 0 | 0 | 81 | 42 | 105 | 75 | 310 | 0 | 0 | 2445 | 15.60 | 39547 |
| 16 | 455 | 455 | 62.37 | 0 | 0 | 0 | 0 | 0 | 81 | 43 | 105 | 75 | 291 | 0 | 0 | 7274 | 46.39 | 37013 |
| 17 | 455 | 410.36 | 20 | 0 | 0 | 0 | 0 | 0 | 81 | 48 | 105 | 75 | 288 | 0 | 0 | 2253 | 14.35 | 35625 |
| 18 | 455 | 455 | 86.35 | 0 | 0 | 0 | 0 | 0 | 82 | 49 | 104 | 75 | 299 | 0 | 0 | 4290 | 27.36 | 37919 |
| 19 | 455 | 455 | 130 | 28.43 | 0 | 0 | 0 | 0 | 82 | 46 | 104 | 75 | 311 | 1815 | 11.57 | 0 | 0 | 39217 |
| 20 | 455 | 455 | 130 | 130 | 55 | 42.7 | 0 | 0 | 82 | 48 | 105 | 75 | 326 | 361 | 2.31 | 0 | 0 | 43911 |
| 21 | 455 | 455 | 130 | 65.18 | 20 | 0 | 0 | 0 | 82 | 47 | 104 | 75 | 317 | 7187 | 45.83 | 0 | 0 | 40612 |
| 22 | 455 | 450 | 0 | 0 | 47.52 | 0 | 0 | 0 | 81 | 43 | 105 | 75 | 297 | 3683 | 23.49 | 0 | 0 | 37159 |
| 23 | 455 | 285.35 | 0 | 0 | 0 | 0 | 0 | 0 | 81 | 44 | 105 | 75 | 278 | 5435 | 34.65 | 0 | 0 | 32793 |
| 24 | 455 | 195.4 | 0 | 0 | 0 | 0 | 0 | 0 | 81 | 40 | 105 | 75 | 262 | 4486 | 28.61 | 0 | 0 | 30704 |